\documentclass{amsart}

\usepackage{amsmath,amsthm,amsfonts,amscd,amssymb}
\usepackage[pdfborder={0 0 0}]{hyperref}
\usepackage{enumerate,verbatim}
\usepackage{enumitem}
\usepackage{color, ulem}

\newcommand{\ZZ}{\mathbb{Z}}

\newtheorem{thm}{Theorem}

\newtheorem{lem}[thm]{Lemma}
\newtheorem{prop}[thm]{Proposition}

\theoremstyle{definition}

\theoremstyle{remark}

\theoremstyle{remark}

\begin{document}

\title{Lattices and quadratic forms from tight frames in Euclidean spaces}
\author{Albrecht B\"ottcher}
\author{Lenny Fukshansky}\thanks{Fukshansky acknowledges support by the NSA grant H98230-1510051 and Simons Foundation grant \#519058}

\address{Fakult\"at f\"ur Mathematik, TU Chemnitz, 09107 Chemnitz, Germany}
\email{aboettch@mathematik.tu-chemnitz.de}
\address{Department of Mathematics, 850 Columbia Avenue, Claremont McKenna College, Claremont, CA 91711, USA}
\email{lenny@cmc.edu}


\begin{abstract}
This paper supplies additions to our paper in Linear Algebra Appl. 510 (2016) 395--420 on integral spans of tight frames in Euclidean spaces. In that previous paper, we considered the case of an equiangular tight frame (ETF), proving that if its integral span is a lattice then the frame must be rational, but overlooking a simple argument in the reverse direction. Thus our first result here is that the integral span of an ETF is a lattice if and only if the frame is rational. Further, we discuss conditions under which such lattices are eutactic and perfect and, consequently, are  local maxima of the packing density function in the dimension of their span. In particular, the unit (276, 23) equiangular tight frame is shown to be eutactic and perfect. More general tight frames and their norm-forms are considered as well, and definitive results are obtained in dimensions two and three.
\end{abstract}

\maketitle

\def\A{{\mathcal A}}
\def\AA{{\mathfrak A}}
\def\B{{\mathcal B}}
\def\C{{\mathcal C}}
\def\D{{\mathcal D}}
\def\EE{{\mathfrak E}}
\def\F{{\mathcal F}}
\def\G{{\mathcal G}}
\def\x{{\mathcal H}}
\def\I{{\mathcal I}}
\def\II{{\mathfrak I}}
\def\J{{\mathcal J}}
\def\K{{\mathcal K}}
\def\kk{{\mathfrak K}}
\def\L{{\mathcal L}}
\def\LL{{\mathfrak L}}
\def\M{{\mathcal M}}
\def\mm{{\mathfrak m}}
\def\MM{{\mathfrak M}}
\def\N{{\mathcal N}}
\def\O{{\mathcal O}}
\def\OO{{\mathfrak O}}
\def\PP{{\mathfrak P}}
\def\R{{\mathcal R}}
\def\W{{\mathcal W}}
\def\PNR{{\mathcal P_N(\real)}}
\def\PMNR{{\mathcal P^M_N(\real)}}
\def\PdNR{{\mathcal P^d_N(\real)}}
\def\s{{\mathcal S}}
\def\V{{\mathcal V}}
\def\X{{\mathcal X}}
\def\Y{{\mathcal Y}}
\def\Z{{\mathcal Z}}
\def\Z{{\mathcal Z}}
\def\({{\langle}}
\def\){{\rangle}}
\def\cee{{\mathbb C}}
\def\Nn{{\mathbb N}}
\def\pee{{\mathbb P}}
\def\que{{\mathbb Q}}
\def\QQ{{\mathbb Q}}
\def\real{{\mathbb R}}
\def\RR{{\mathbb R}}
\def\zed{{\mathbb Z}}
\def\ZZ{{\mathbb Z}}
\def\aaa{{\mathbb A}}
\def\ff{{\mathbb F}}
\def\HDelta{{\it \Delta}}
\def\kk{{\mathfrak K}}
\def\qbar{{\overline{\mathbb Q}}}
\def\kbar{{\overline{K}}}
\def\ybar{{\overline{Y}}}
\def\kkbar{{\overline{\mathfrak K}}}
\def\ubar{{\overline{U}}}
\def\eps{{\varepsilon}}
\def\ahat{{\hat \alpha}}
\def\bhat{{\hat \beta}}
\def\k{{\nu}}
\def\gt{{\tilde \gamma}}
\def\h{{\tfrac12}}
\def\be{{\boldsymbol e}}
\def\bei{{\boldsymbol e_i}}
\def\bc{{\boldsymbol c}}
\def\bdt{{\boldsymbol \delta}}
\def\bff{{\boldsymbol f}}
\def\bm{{\boldsymbol m}}
\def\bk{{\boldsymbol k}}
\def\bi{{\boldsymbol i}}
\def\bl{{\boldsymbol l}}
\def\bq{{\boldsymbol q}}
\def\bu{{\boldsymbol u}}
\def\bt{{\boldsymbol t}}
\def\bs{{\boldsymbol s}}
\def\bv{{\boldsymbol v}}
\def\bw{{\boldsymbol w}}
\def\bx{{\boldsymbol x}}
\def\bbx{{\overline{\boldsymbol x}}}
\def\bX{{\boldsymbol X}}
\def\bz{{\boldsymbol z}}
\def\bwy{{\boldsymbol y}}
\def\bY{{\boldsymbol Y}}
\def\bL{{\boldsymbol L}}
\def\ba{{\boldsymbol a}}
\def\bb{{\boldsymbol b}}
\def\bet{{\boldsymbol\eta}}
\def\bxi{{\boldsymbol\xi}}
\def\bo{{\boldsymbol 0}}
\def\bone{{\boldsymbol 1}}
\def\bol{{\boldsymbol 1}_L}
\def\ep{\varepsilon}
\def\p{\boldsymbol\varphi}
\def\q{\boldsymbol\psi}
\def\rank{\operatorname{rank}}
\def\aut{\operatorname{Aut}}
\def\lcm{\operatorname{lcm}}
\def\sgn{\operatorname{sgn}}
\def\spn{\operatorname{span}}
\def\md{\operatorname{mod}}
\def\Norm{\operatorname{Norm}}
\def\dim{\operatorname{dim}}
\def\det{\operatorname{det}}
\def\Vol{\operatorname{Vol}}
\def\rk{\operatorname{rk}}
\def\ord{\operatorname{ord}}
\def\ker{\operatorname{ker}}
\def\div{\operatorname{div}}
\def\Gal{\operatorname{Gal}}
\def\GL{\operatorname{GL}}
\def\SNR{\operatorname{SNR}}
\def\WR{\operatorname{WR}}
\def\IWR{\operatorname{IWR}}
\def\scg{\operatorname{\left< \Gamma \right>}}
\def\swrh{\operatorname{Sim_{WR}(\Lambda_h)}}
\def\ch{\operatorname{C_h}}
\def\cht{\operatorname{C_h(\theta)}}
\def\scgt{\operatorname{\left< \Gamma_{\theta} \right>}}
\def\scgmn{\operatorname{\left< \Gamma_{m,n} \right>}}
\def\gat{\operatorname{\Omega_{\theta}}}
\def\Obar{\operatorname{\overline{\Omega}}}
\def\Lbar{\operatorname{\overline{\Lambda}}}
\def\mn{\operatorname{mn}}
\def\disc{\operatorname{disc}}
\def\rot{\operatorname{rot}}
\def\Prob{\operatorname{Prob}}
\def\co{\operatorname{co}}
\def\ot{\operatorname{o_{\tau}}}
\def\Aut{\operatorname{Aut}}
\def\Mat{\operatorname{Mat}}
\def\SL{\operatorname{SL}}
\def\id{\operatorname{id}}

\section{Introduction and Main Results}
\label{intro}

\noindent
We denote by $\left<\ ,\ \right>$ the usual inner product on $\real^k$  and by $\|\bx\| := \left< \bx,\bx \right>^{1/2}$ the Euclidean norm.
Let $n \ge k$ and let $\F := \left\{ \bff_1,\dots,\bff_n \right\} \subset \real^k$ be a set of vectors such that
$ \spn_{\real} \left\{ \bff_1,\dots,\bff_n \right\} =\real^k$.
Put
$$\Lambda(\F) = \spn_{\zed} \left\{ \bff_1,\dots,\bff_n \right\}.$$
We consider $\bff_1,\dots,\bff_n$ as column vectors and denote by $G$ the $k \times n$ matrix with these
vectors as columns,
Clearly, we may think of $\Lambda(\F)$ as the set $\{G\ba: \ba \in \zed^n\}$.
The {\it norm-form} associated with $\F$
is the quadratic form
\[Q_\F(\ba)=\| G\ba\|^2=\langle G^\top G \ba, \ba \rangle.\]
A quadratic form $Q(\ba)=\langle H \ba, \ba \rangle$ with a symmetric matrix $H$
is said to {\it have separated values} if the values of $Q(\ba)$ for  $\ba$ in $\zed^n$ are separated by a positive
number, that is, if
\[\inf\{|Q(\ba)-Q(\bb)|: \ba,\bb \in \zed^n, Q(\ba) \neq Q(\bb)\} > 0.\]
We call the set $\F$ {\it rational} if the inner products $\left< \bff_i, \bff_j \right>$ are rational numbers for all $1 \leq i,j \leq n$.
This is equivalent to saying that the entries of the $n \times n$ Gram matrix $G^\top G$ are all rational.

We here study two questions. First, we are interested in conditions ensuring that $\Lambda(\F)$ is a lattice, and if it is, in properties of this lattice.
Secondly, we look for conditions ensuring that $Q_\F$ has separated values. Recall that a {\it lattice} is a discrete additive subgroup of $\real^k$.
The set $\Lambda(\F)$ of the integer linear combinations of the vectors $\bff_1, \ldots, \bff_n$ is clearly an additive subgroup of $\real^k$,
but it may contain accumulation points and hence not be discrete, which would prevent it from being a lattice. Obviously, $\Lambda(\F)$ is
a lattice if and only if $0$ is not an accumulation point of the values of $Q_\F$, that is, if and only if there exists an $\varepsilon >0$ such that
$Q_\F(\ba) \notin (0,\varepsilon)$ for $\ba \in \zed^n$. A lattice $\Lambda \subset \real^k$ is said to be of full rank if $\spn_\real \Lambda =\real^k$.
As we always suppose that
$\spn_{\real} \{ \bff_1,\dots,\bff_n \}=\real^k$, the set $\Lambda(\F)$ is a lattice if and only if it is
a full-rank lattice.

For further reference, we state the following simple observation, which, unfortunately, was overlooked
in~\cite{etf1}.

\begin{prop} \label{Thm11}
If $\F$ is rational, then $\Lambda(\F)$ is a lattice and the values of $Q_\F$ are separated.
\end{prop}

\proof There is an integer $d >0$ such that all entries of $dG^\top G$ are integers, and hence $dQ_\F(\ba)=d\langle G^\top G\ba,\ba\rangle$
assumes values in $\{0,1,2,\ldots\}$ for $\ba \in \zed^n$, which shows that $Q_\F$ has separated values and does not take values in $(0,1/d)$.
\endproof

Clearly, if $\mu \neq 0$ is any real number, then
$\Lambda(\mu\F$) is also a lattice and the values of $Q_{\mu \F}$ are separated.

Much is known in the case $n=k$. Suppose $G$ is an invertible square matrix. This ensures that $\Lambda(\F)$
is a lattice. We say that $\F$ is {\it irrational} if $G^\top G$ is not a non-zero scalar multiple of
a matrix with rational entries. Note that an {\it irrational} $\F$ is not the exact logical negation of
a {\it rational} $\F$. In fact, we have exactly three mutually excluding possibilities. First, $\F$ is {\it rational}
if all entries of $G^\top G$ are rational. Otherwise $\F$ is called {\it non-rational}, and in that case we have exactly two
possibilities:
either there is a non-zero scalar $\mu$ such $\mu^2 G^\top G$ has rational entries or $\mu^2 G^\top G$ has an irrational
entry for each non-zero scalar $\mu$. In the last case, $\F$ is called {\it irrational}.

An easy example of an irrational set $\F$ for which $Q_\F$ does not have separated
values is given by the $3 \times 3$ matrix
\[G=\left(\begin{array}{lll} 1 & \xi & 0\\
0 & 1 & 0\\
0 & 0 & 1 \end{array} \right),\]
where $\xi$ is irrational. In that case
\[Q_\F(\ba)=(a_1+\xi a_2)^2+ a_2^2+a_3^2.\]
Given any $\varepsilon >0$, there are integers $a_1,a_2$ such that $0<(a_1+\xi a_2)^2 < \varepsilon$.
Thus, $(a_1+\xi a_2)^2+ a_2^2=a_2^2+\delta$ with $\delta \in (0,\varepsilon)$.
It follows that
\[Q_\F(a_1,a_2,0)-Q_\F(0,0,a_2)=\delta \in (0,\varepsilon),\]
which shows that the values of $Q_\F$ are not separated.
A famous result by Margulis~\cite{marg,margulis} implies that if $\F$ is irrational and $k \ge 2$, then the values
of $Q_\F(\ba)-Q_\F(\bb)$ for $\ba,\bb \in \zed^k$ are dense in $\real$, implying that the
values of $Q_\F$ are not separated.

In the case $n >k$ it may also happen that $\Lambda(\F)$ is a lattice and the values of the associated norm form
$Q_\F$ are not separated. Indeed, let $G=(G_0\;G_1)$ with an invertible $k \times k$ matrix and a $k \times (n-k)$
matrix $G_1$. Suppose $G_0^\top G_0$ is irrational and each column of $G_1$ is an integer linear combination of the columns of $G_0$. Then $\Lambda(G)$
is a lattice (with $G_0$ as a basis matrix), but, again by the result of Margulis, the values of indefinite irrational
quadratic form $Q_\F(\ba)-Q_\F(\bb)$ with $\ba,\bb \in \zed^n$  are dense in $\real$, which implies that the values of $Q$
are not separated.

In the following, we treat the case $n >k$ and special sets $\F$, so-called tight frames.
The set $\F = \left\{ \bff_1,\dots,\bff_n \right\} \subset \real^k$
is called an $(n,k)$ {\it tight frame} if $n > k$, $\spn_{\real} \{ \bff_1,\dots,\bff_n \}=\real^k$
and there exists a positive number $\gamma \in \real$ such that
\begin{equation}
\|\bx\|^2 =  \gamma \sum_{i=1}^n \left< \bff_i, \bx \right>^2\label{gamma}
\end{equation}
for every $\bx \in \real^k$.  Obviously, requirement~(\ref{gamma}) is equivalent to the equality $G G^\top = \gamma I$, which in turn simply means
that the $k$ rows of $\gamma^{-1/2}G$ are rows of an orthogonal $n \times n$ matrix.

Clearly, if $\F$ is a rational $(n,k)$ tight frame and $\mu \neq 0$ is a real number such that $\mu^2 \notin \que$,
then $\mu\F=\{\mu \bff_1, \ldots, \mu \bff_n\}$ is no longer rational but still an $(n,k)$ tight frame for which $\Lambda(\mu \F)$ is a lattice
and for which $Q_\F$ has separated values. Except for this trivial construction of non-rational tight frames generating a lattice,
we do not know any non-rational tight frames that induce a lattice or a norm-form with separated values. We are able to prove
the following theorem, which settles the case of tight frames in $\real^2$ and $\real^3$. Note that requiring that one of the vectors $\bff_1, \ldots, \bff_n$ has length $1$ rules out multiplication by $\mu$ with $\mu^2 \notin \que$.

\begin{thm} \label{Thm12}
Let $\F=\{\bff_1,\ldots,\bff_n\}$ be an $(n,2)$ or an $(n,3)$ tight frame containing at least one unit vector. Then
the following are equivalent:

\smallskip
(i) $\F$ is rational,

\smallskip
(ii) $\Lambda(\F)$ is a lattice,

\smallskip
(iii) $Q_\F$ has separated values.
\end{thm}

Note that the implication $(iii) \Rightarrow (ii)$ is trivial, because $0$ cannot be an accumulation point of
a separated set. The implications $(i) \Rightarrow (ii)$ and $(i) \Rightarrow (iii)$ follow from Proposition~\ref{Thm11}. Thus, we are left with the implication $(ii) \Rightarrow (i)$, which will be proved in this paper.

A {\it unit $(n,k)$ equiangular tight frame} (ETF) of angle $c \in (0,1)$ is an $(n,k)$ tight frame $\F=\{\bff_1, \ldots,\bff_n\}$
consisting of unit vectors such that $\left| \left< \bff_i, \bff_j \right> \right| = c$ for every pair $1 \leq i \neq j \leq n$.
It is well known that if $\F$ is a unit $(n,k)$ ETF, then the constant $\gamma$ in~(\ref{gamma}) necessarily equals $k/n$, that is,
we have
$$\|\bx\|^2 =  \frac{k}{n} \sum_{i=1}^n \left< \bff_i, \bx \right>^2$$
for every $\bx \in \real^k$.
We can write $c=1/\alpha$ with $\alpha \in [1, \infty)$, and it is known, e.g., from~\cite{Sus}, that
\begin{equation}
c = \frac{1}{\alpha} = \sqrt{ \frac{n-k}{k(n-1)}}. \label{calph}
\end{equation}
Moreover, in~\cite{Sus} it is shown that if $n \neq k+1$ or $n \neq 2k$, then $\alpha$ is an odd integer.
On the other hand, if $n=k+1$ or $n=2k$, then~(\ref{calph}) implies that $\alpha$ is either an integer or irrational.
The Gram matrix $G^\top G$
of a unit $(n,k)$ ETF $\F$ has $1$ on the main diagonal and $\pm 1/\alpha$ elsewhere. Thus, $\F$
is a rational frame if and only if $\alpha$ is a positive integer.
We here prove the following.

\begin{thm} \label{Thm13} Let $\F$ be a unit $(n,k)$ ETF. Then the following are equivalent:

\smallskip
(i) $\alpha$ is a positive integer,

\smallskip
(ii) $\Lambda(\F)$ is a lattice,

\smallskip
(iii) $Q_\F$ has separated values.

\noindent
If $\alpha$ is a positive integer, then the minimal norm $|\Lambda(\F)|$ of $\Lambda(\F)$ satisfies
\begin{equation}
\label{min-norm}
|\Lambda(\F)| := \min \{ \|\bx\| : \bo \neq \bx \in \Lambda(\F) \} \geq \frac{1}{\sqrt{\alpha}}.
\end{equation}
\end{thm}

Notice that $k$ in Theorem~\ref{Thm13} may be arbitrary, but we require equiangularity, whereas in Theorem~\ref{Thm12} the value of $k$ is restricted to $2$ and $3$, but no equiangularity is needed.
The implication $(ii) \Rightarrow (i)$ in Theorem~\ref{Thm13} was proved in~\cite{etf1}. Since, as said, a unit ETF is a rational frame if and only if $\alpha$ is a rational number,
Proposition~\ref{Thm11} yields the implications $(i) \Rightarrow (ii)$ and $(i) \Rightarrow (iii)$ in Theorem~\ref{Thm13}.
The implication $(iii) \Rightarrow (ii)$ of this theorem is trivial, and hence we are left with proving the estimate~(\ref{min-norm}).

Since it is known that $\alpha$ is rational for  $n \neq 2k$, it follows that all unit $(n,k)$ ETFs with $n \neq 2k$ generate lattices with separated values
in~$\real^k$. There are also infinitely many unit $(n,k)$ ETFs with $n = 2k$ and $\alpha$ rational; see~\cite{Goet} and~\cite{Sus}.

Let
$$S(\Lambda(\F)) := \{ \bx \in \Lambda(\F) : \|\bx\| = |\Lambda(\F)| \}$$
be the set of minimal vectors of the lattice $\Lambda(\F)$.
A collection of points $X = \{ \bx_1,\dots,\bx_m \}$ on the unit sphere in $\real^k$, $m \geq k$, is called {\it weakly eutactic} if there exist real numbers $c_1,\dots,c_m$, called {\it eutaxy coefficients}, such that
$$\|\bwy\|^2 = \sum_{i=1}^m c_i \left< \bwy,\bx_i \right>^2$$
for every vector $\bwy \in \real^k$. The set $X$ is called {\it semi-eutactic} if $c_1,\dots,c_m \geq 0$, {\it eutactic} if $c_1,\dots,c_m > 0$, and {\it strongly eutactic} if $c_1 = \dots = c_m > 0$. It is easy to observe that if $X$ is a set containing a weakly eutactic or semi-eutactic subset $Y$, then $X$ itself is weakly or semi-eutactic, respectively.
These notions from lattice theory are closely related to the notions of scalability used in frame theory. For example,
in the terminology of~\cite{Kut}, the semi-eutactic sets $X$ are the scalable sets and the eutactic sets are the strictly scalable ones.
Notice in particular that a unit $(n,k)$ ETF $\F$ is a strongly eutactic set with the eutaxy coefficients equal to $k/n$, and the set $\pm \F$ is strongly eutactic with the eutaxy coefficients equal to $k/2n$. A lattice is called weakly eutactic, semi-eutactic, eutactic, or strongly eutactic if its set of minimal vectors has the respective property. Based on a number of examples we have worked out, we conjecture that $|\Lambda(\F)| = 1$ and $S(\Lambda(\F)) = \pm \F$ for every unit ETF. If this is the case, it immediately follows that the lattices $\Lambda(\F)$ strongly eutactic; see Proposition~2.4 of~\cite{etf1}.

Further, a collection of points $X = \{ \bx_1,\dots,\bx_m \}$ (written as column-vectors) on the unit sphere in $\real^k$, $m \geq k$, is called {\it perfect} if the set of real symmetric matrices
$$X^* := \left\{ \bx_1 \bx_1^\top, \dots, \bx_m \bx_m^\top \right\}$$
spans the entire space of real symmetric $k \times k$ matrices as a real vector space. It is easy to see that if $X$ is a set containing a perfect subset $Y$, then $X$ itself is perfect. In fact, Theorem~3.6.2 on p.~85 of~\cite{martinet} states that if $X$ contains a perfect and eutactic subset $Y$, then $X$ itself is perfect and eutactic. A lattice $\Lambda$ is called perfect if its set of minimal vectors $S(\Lambda)$ is perfect. A classical theorem of Voronoi (1908) states that a lattice is perfect and eutactic if and only if it is extreme, i.e., is a local maximum of the packing density function in its dimension; see Theorem~3.4.6 on p.~81 of~\cite{martinet}.

The well-known Gerzon bound for a unit $(n,k)$ ETF states that
\begin{equation}
\label{eq1}
k < n \leq \frac{k(k+1)}{2}.
\end{equation}
We will say that a unit $(n,k)$ ETF is {\it maximal} if Gerzon's upper bound is attained. While unit $(n,k)$ ETFs exist in every dimension,
there are only four known examples of maximal ETFs: in dimensions 2, 3, 7, and 23. We here prove the following.

\begin{thm} \label{Thm14} Let $k > 3$ and $\F$ be a maximal unit $(n,k)$ ETF. Assuming $|\Lambda(\F)| = 1$, the lattice $\Lambda(\F)$ is perfect and eutactic, and hence extreme.
\end{thm}

As just mentioned, there are only four maximal unit ETFs currently known: these are $(3,2)$, $(6,3)$, $(28,7)$, and $(276,23)$ ETFs.
The $(3,2)$ ETF generates the hexagonal lattice, whose minimal vectors are $\pm$ the frame vectors, and hence it is perfect and strongly eutactic (and has, independently
of these two properties, been well-known to be extreme since Gauss). The $(6,3)$ ETF has irrational $\alpha$, and hence, as already shown in ~\cite{etf1}, does not generate a lattice.
Also in~\cite{etf1}, we proved that
the $(28,7)$ ETF generates a lattice with $\pm$ the frame vectors as its minimal vectors and then showed in a very computational manner
that this lattice is perfect and strongly eutactic and thus extreme. Theorem~\ref{Thm14} now gives the conclusion straight away.
The remaining case is that of the lattice generated by the $(276,23)$ ETF. Here is our result.

\begin{thm} \label{Thm15} Let $\F$ be the unit $(276,23)$ ETF. Then $\Lambda(\F)$ is perfect and eutactic, and hence extreme.
\end{thm}

We present the proofs of the above results in Section~\ref{proofs}.
\bigskip

\section{Proofs}
\label{proofs}

\medskip
\proof[Proof of Theorem~\ref{Thm12}]
As said, we are left with the implication $(ii) \Rightarrow (i)$.

We first consider the case $k=2$. Suppose $\Lambda(\F)$ is a lattice.
Orthogonal linear maps transform tight frames into tight frames and lattices into lattices,
and they preserve inner products. We may therefore
without loss of generality assume that $\bff_1=(1,0)^\top.$ Thus, the $2 \times n$ matrix $G$
may be written as
\[G=\left(\begin{array}{ccccc} 1 & a & x_1 & \ldots & x_{n-2}\\0 & b & y_1 & \ldots & y_{n-2}\end{array}\right).\]
We may assume that $b \neq 0$. Condition (\ref{gamma}) is equivalent to the two equalities
\begin{equation}
1+a^2+\sum x_i^2=b^2+\sum y_i^2, \qquad ab+\sum x_iy_i=0, \label{two}
\end{equation}
all sums being from $i=1$ to $i=n-2$.
Proposition 3.1 of \cite{etf1} implies that $\Lambda(\F)$ is a lattice if and only if the $n-2$ vectors
\[\left(\begin{array}{cc} 1 & a\\ 0 & b\end{array}\right)^{-1}\left(\begin{array}{c} x_i\\ y_i \end{array}\right)
=\left(\begin{array}{c}x_i-a \frac{y_i}{b} \\[1ex] \frac{y_i}{b}\end{array}\right)\]
are rational: while this proposition is stated for unit tight frames, the proof shows that it is also valid
for arbitrary tight frames. It follows that $q_i:=y_i/b$ are rational numbers and that we have
$x_i=aq_i+p_i$ with rational numbers $p_i$. Inserting this into the second equality of~(\ref{two})
we get
\[a\left(1+\sum q_i^2\right)+\sum p_iq_i =0,\]
which shows that $a$ is rational. The first equality of~(\ref{two}) gives
\[1+a^2+\sum (aq_i+p_i)^2=b^2\left(1+\sum q_i^2\right),\]
and this implies that $b^2$ is rational. Using that $a, q_i,p_i,b^2$ are rational numbers, it is
easily verified that all inner products $\langle \bff_j, \bff_k \rangle$ are rational numbers,
which means that $\F$ is a rational frame.

Now suppose $k=3$ and $\Lambda(\F)$ is a lattice. Again using an orthogonal transformation,
we may assume that the first three columns of $G$ form an upper-triangular invertible matrix
and that the first column is $(1,0,0)^\top$. Let $(c,\mu,0)^\top$ be the second column.
Since the upper-triangular matrix is invertible, we have $\mu \neq 0$ and may therefore write
the entire matrix $G$ in the form
\[G=\left(\begin{array}{cccccc}
1 & c & d & z_1 & \ldots & z_{n-3}\\
0 & \mu & \mu a & \mu x_1 & \ldots & \mu x_{n-3}\\
0 & 0 & \mu b & \mu y_1 & \ldots & \mu y_{n-3} \end{array}\right).\]
From Proposition 3.1 of \cite{etf1} we deduce that the $n-3$ columns
\[\left(\begin{array}{ccc} 1 & c & d\\
0 & \mu & \mu a\\
0 & 0 & \mu b\end{array}\right)^{-1}\left(\begin{array}{c}
z_i \\ \mu x_i \\ \mu y_i \end{array}\right)
= \left(\begin{array}{c} z_i-cx_i - \left( \frac{d}{b} - \frac{ac}{b} \right)y_i\\[1ex]
x_i-\frac{a}{b}y_i\\[1ex]
\frac{y_i}{b}\end{array}\right) = \begin{pmatrix} r_i \\ p_i \\ q_i \end{pmatrix} \]
are all rational. It follows that
\begin{eqnarray*}
& & y_i=bq_i, \quad  x_i = \frac{a}{b}  b q_i + p_i =aq_i+p_i, \\
& & z_i=c(aq_i+p_i)+\left( \frac{d}{b} - \frac{ac}{b} \right)bq_i+r_i = cp_i +dq_i+r_i,
\end{eqnarray*}
and since~(\ref{two}) holds also in the case at hand, we obtain as above that $a$ and $b^2$
must be rational. Taking the inner product of the first and third rows of $G$ we get
\begin{eqnarray*}
0 & = & \mu bd+\mu \sum y_i z_i\\
& = &\mu bd\left(1+\sum q_i^2\right) +\mu bc\sum p_iq_i +\mu b \sum q_i r_i.
\end{eqnarray*}
Cancelling out $\mu b$, we arrive at the equality
\[0 = d\left(1+\sum q_i^2\right)+c \sum p_iq_i +\sum q_i r_i,\]
which shows that $d =cs_1+s_2$ with rational
$$s_1 = -\frac{\sum q_i p_i}{1 + \sum q_i^2},\ s_2 = - \frac{\sum q_i r_i}{1 + \sum q_i^2}.$$
This gives $z_i=ct_i+u_i$ with rational numbers $t_i = p_i+s_1q_i$ and $u_i = r_i+s_2q_i$. Consequently,
\begin{eqnarray}
\label{s1+}
s_1 + \sum t_i q_i & = & s_1 + \sum (p_i+s_1q_i) q_i = s_1 \left( 1 + \sum q_i^2 \right) + \sum p_i q_i \nonumber \\
& = & \left( -\frac{\sum q_i p_i}{1 + \sum q_i^2} \right) \left( 1 + \sum q_i^2\right) + \sum p_i q_i = 0.
\end{eqnarray}
Thus, at this point we have
\begin{equation}
G=\left(\begin{array}{cccccc}
1 & c & cs_1+s_2 & \ldots & ct_i+u_i & \ldots\\
0 & \mu & \mu a & \ldots & \mu(aq_i+p_i) & \ldots \\
0 & 0 & \mu b & \ldots & \mu b q_i & \ldots \end{array}\right).\label{Gfin}
\end{equation}
After cancelling out $\mu b$, the inner product of the first and third rows of $G$ shows that
\begin{equation}
0 =  c \left(s_1+\sum t_iq_i\right) + s_2 + \sum u_i q_i = s_2 + \sum u_i q_i,\label{E2}
\end{equation}
by~\eqref{s1+}. After cancelling out $\mu$, the inner product of the first and second rows of $G$ yields
\begin{eqnarray*}
0 & = & c + a(cs_1+s_2) + \sum (ct_i+u_i)(aq_i+r_i) \\
& = & c \left( 1 + a (s_1+\sum t_i q_i) + \sum t_i r_i \right) + a \left( s_2 + \sum u_i q_i \right) + \sum u_i r_i 
\end{eqnarray*}
and since $s_1+\sum t_iq_i=0$ and $s_2+\sum u_iq_i=0$ by~\eqref{s1+} and~\eqref{E2}, it follows that
\begin{equation}
0= c \left( 1 + \sum t_i r_i \right) + \sum u_i r_i. \label{E1}
\end{equation}
We claim that $1 + \sum t_i r_i \neq 0$. To see this, assume the contrary, that is, let $1 + \sum t_i r_i = 0$. Then~\eqref{E1} implies that $\sum u_i r_i = 0$. Recalling that $r_i = u_i - s_2q_i$ and $\sum u_i q_i = -s_2$, we obtain:
$$0 = \sum u_i (u_i - s_2q_i) = \sum u_i^2 - s_2 \sum u_i q_i = \sum u_i^2 + s_2^2,$$
which is only possible if $s_2 = 0$ and $u_i=0$ for all $i$. This means that $r_i = u_i - s_2q_i = 0$ for all $i$, and hence
$$1 + \sum t_i r_i = 1,$$
which contradicts our assumption. Thus, $1 + \sum t_i r_i$ is indeed nonzero. So~\eqref{E1} implies that $c$ is rational.
Finally, since the first and third rows of $G$ have equal norm,
we obtain
\[1+c^2+(cs_1+s_2)^2 +\sum(ct_i+u_i)^2=\mu^2 b^2 \left(1+\sum q_i^2\right),\]
which tells us that $\mu^2$ is rational.

In summary, $G$ is of the form~(\ref{Gfin}) with $a, c, b^2, \mu^2,
s_1, s_2, t_i, u_i, q_i, r_i$ rational. Taking this into account it can be
checked directly that the inner products of any two columns of $G$ are rational
numbers. Hence $\F$ is rational.
\endproof

\proof[Proof of Theorem~\ref{Thm13}]  It suffices to prove~\eqref{min-norm}.
We have $\Lambda(\F) = \{ G \ba : \ba \in \zed^n \}$. The $j,k$ entry of $G^\top G$ is $\langle \bff_j,\bff_k\rangle=\pm 1/\alpha$
for $j \neq k$ and $1$ for $j=k$.  Consequently, $\alpha G^\top G$ is an integer matrix. It follows
that $\alpha\|G\ba\|^2 = \alpha \left< G^{\top} G \ba, \ba \right>$ has its values in $\{0,1,2,\ldots\}$. Thus,
$\| G\ba\|^2 \ge 1/\alpha$ whenever $G\ba \neq \bo$.
\endproof

\begin{lem} \label{perfect_ETF} Let $\F = \{ \bff_1,\dots,\bff_n \}$ be a maximal unit $(n,k)$ ETF in $\real^k$, that is, $n = k(k+1)/2$.
Then $\F$ is a perfect set in $\real^k$.
\end{lem}

\proof
We only need to show that the set
$$\F^* := \left\{ \bff_1 \bff_1^\top, \dots, \bff_n \bff_n^\top \right\}$$
is linearly independent. The following argument comes from the proof of Gerzon's bound, as outlined, for instance, at the top of the second page of~\cite{benny}.
Treating the matrices $\bff_i \bff_i^\top$ as vectors in $\real^{n^2}$, one can compute the inner products as
$$\left< \bff_i \bff_i^\top, \bff_j \bff_j^\top \right> = \left< \bff_i, \bff_j \right>^2,$$
which is $1/\alpha^2$ when $i \neq j$ and 1 when $i = j$.
Thus, the Gram matrix of $\F^*$ is
\[\left(1-\frac{1}{\alpha^2}\right) I+\frac{1}{\alpha^2} J,\]
where $J$ is the $n \times n$ matrix all entries of which are $1$. The eigenvalues of this Gram matrix
are $1-\alpha^2>0$ and $1-1/\alpha^2+n >0$, so the matrix is invertible, which implies that $\F^*$ is linearly independent.
\endproof

\proof [Proof of Theorem~\ref{Thm14}]
Since $k > 3$ and $\F$ is maximal, we have $n = k(k+1)/2 > 2k$. Hence $\alpha \in \zed$, and so $\Lambda(\F)$ is a lattice
by Proposition~\ref{Thm11}. Assuming $|\Lambda(\F)| = 1$, we see that $\pm \F \subset S(\Lambda(F))$. This inclusion implies that $\Lambda(\F)$ is strongly eutactic
and Lemma~\ref{perfect_ETF} tells us that the lattice is perfect.
\endproof

\proof [Proof of Theorem~\ref{Thm15}]
Let $\L$ be the Leech lattice. There is a normalization constant $\mu >0$ such that
the vectors of the unit $(276,23)$ ETF $\F$ appear among the minimal vectors of the normalized Leech lattice $\mu\L$;
see Section~16.3 of~\cite{deza} as well as a reference to it in Section~3.1 of~\cite{achill1}. Thus, $\F \subset S(\mu \L)$.
It follows that $|\mu\L| \ge 1$ and that $\Lambda(\F)$ is a sublattice of $\mu\L$. As the minimal vectors of a sublattice cannot be shorter than those of
the lattice, we obtain that $|\Lambda(\F)| \ge |\mu\L|$. In summary, $|\Lambda(\F)| \ge 1$, and
since $\F \subset S(\Lambda(\F))$, we arrive at the conclusion that $|\Lambda(\F)| = 1$. The theorem is now immediate from Corollary~\ref{Thm14}.
\endproof

\medskip
{\bf Acknowledgement.} We sincerely thank the referee for several insightful comments and useful suggestions
that improved the original version of this note. We are also greatly indebted to Yuxin (Jessie) Xin for detecting some mistakes in our previous  
proof of Theorem 2. These are corrected in the present version.


\end{document}